\documentclass[12pt]{article}

\usepackage{amsmath}
\usepackage{amsfonts}
\usepackage{amsthm}
\usepackage{graphics}
\usepackage{epsf}

\numberwithin{equation}{section}

\def\Cmpx{{\mathbb{C}}}
\def\Real{{\mathbb{R}}}
\def\Intg{{\mathbb{Z}}}
\def\cnj{\overline}

\def\tfrac#1#2{{\textstyle{\frac{#1}{#2}}}}

\def\nfrac#1#2{{\raisebox{.5ex}{$#1$}\!/\!\raisebox{-.5ex}{$#2$}}}
\def\scrhalf{{\raisebox{.3ex}%
{$\scriptstyle 1$}\!/\!\raisebox{-.3ex}{$\scriptstyle 2$}}}

\def\wh{\widehat}

\def\innerprod(#1){{\langle #1 \rangle}}
\def\norm#1{\|#1\|}
\def\rvec(#1){{|#1\rangle}} 
\def\lvec(#1){{\langle#1|}} 

\def\ppmatrix#1{\begin{pmatrix} #1 \end{pmatrix}}

\def\defn#1{{\em #1}}

\def\Weil{{{\cal W}}}
\def\spinset{{{\cal S}}}
\def\PsiSpace{{\Psi}}
\def\Pexpr{{\cal P}}

\def\baseterm(#1,#2,#3){{\psi{}^{(#2,#3)}_{#1}}}

\def\PsirSpace#1{{{\Psi}{}^{(#1,\bullet)}_{\bullet}}}

\def\Rp{{\hat R}}
\def\Ad#1{{\mathrm{ad}}_{#1}}

\def\Wigner(#1,#2,#3){{\left<
\begin{array}{ccc} & #2 & \\ #1 & & 0 \\ & #3 & \end{array}\right>}}
\def\Czz{{C_{00}}}
\def\RMnrk(#1,#2,#3)%
{{\langle #3+#2 \| \baseterm(#1,#2,\bullet) \| #3 \rangle}}
\def\CG(#1,#2,#3,#4,#5,#6){{C^{#1}_{#4}{}^{#2}_{#5}{}^{#3}_{#6}}}
\def\RMpsi(#1,#2,#3,#4,#5,#6){{R^{#1}_{#4}{}^{#2}_{#5}{}^{#3}_{#6}}}
\def\WsixJ(#1,#2,#3,#4,#5,#6)%
{\left\{ \begin{array}{ccc} #1\ &\ #2\ &\ #3
\\#4\ &\ #5\ &\ #6\end{array}\right\}}
\def\Pexpr{{\cal P}}
\def\Wcoeff(#1,#2,#3,#4,#5,#6){{W^{#1}_{#4}{}^{#2}_{#5}{}^{#3}_{#6}}}
\def\Rp{{\wh R}}
\newcounter{bitcnt}
\def\bititem{\refstepcounter{bitcnt}(\arabic{bitcnt})}
\newcounter{lstcnt}
\newenvironment{mylist}{\begin{list}{(\arabic{lstcnt})}%
{\usecounter{lstcnt}\setlength{\rightmargin}{\leftmargin}}}%
{\end{list}}
\def\manif{{\cal M}}

\title{An Introduction to the Noncommutative Sphere and some
Extensions} 

\author{{\bf J. Gratus%
\thanks{Funded by the Royal Society of London European 
Science Exchange Programme}}%
\\
Laboratoire de Gravitation et Cosmologie Relativistes%
\thanks{\it Laboratoire associ\'e au CNRS {\rm URA 769}}
\\
Tour 22/12 4eme etage, Boite Courrier142, 4pl Jussieu. F75252
Paris
\\
email: gratus@ccr.jussieu.fr
}

\date{October 8, 1997}

\begin{document}

\maketitle

\begin{abstract}
This is a copy of the talk given at the conference ``Methods in Field
Theory'' at Star\'a Lesn\'a, The Slovak Republic. Sepemeber 22-26,
1997.  An introduction to the noncommutative sphere and a summary of
the results of articles q-alg/9703038 \cite{Gratus5}, and
q-alg/9708003 \cite{Gratus6} is given. This includes results about the
the algebra of scalar, spinor and vector fields on the noncommutative
sphere.  Possible extensions of these results including a ``Wick
rotation'' to the one and two sheeted hyperboloid are also examined.
\end{abstract}

\tableofcontents

\newpage

\section{Introduction}

Noncommutative geometry is, at present, similar to string theory during
the 80's. If you ask five noncommutative geometers what they mean by
the subject you will get at least half a dozen different answers.
Here, at least for me, is a working definition:

\begin{quote}
Almost all the tools{}\footnote{One can generate the original manifold
from the algebra of functions on that manifold only if the manifold is
algebraic or holomorphic. However the detailed differential structure
of a manifold is seldom of interest to physicists who usually assume
at least a $C^\infty$ structure.} used to study ordinary differential
geometry can be expressed in terms of the functions from a base
manifold to $\Real$ or $\Cmpx$. The algebra of these functions using
pointwise multiplication is commutative.  The object of Noncommutative
Geometry is to replace this algebra with a noncommutative algebra, and
then to find and interpret analogues of these tools.
\end{quote}

The rest of noncommutative geometry is, in my opinion, a wish list.
These include:
\begin{mylist}
\item One would like a one parameter algebra ${\cal A}(\hbar)$, such
that ${\cal A}(\hbar=0)=C({\manif})$, the commutative algebra of
functions on the manifold with pointwise multiplication, and such that
the commutator is of order $\hbar$. I.e. $[f,g]=fg-gf=O(\hbar)$.

\item If $\manif$ has a Poisson structure then one would like the
antisymmetric part of the first order correction to agree with the
Poisson bracket. $\{f,g\}=\lim_{\hbar\to0}(\tfrac1\hbar [f,g])$.
In this sense it is similar to deformation quantisation.

\item A definition of vector fields, and covector fields over $\manif$.

\item A definition of exterior derivative and corresponding exterior
calculus.

\item A definition of intergration over the manifold, and thus a
Hilbert space structure for $C({\manif})$

\item Spinor fields if $\manif$ admits a spinor structure.

\item If $\manif$ is to be thought of as spacetime we need to define
objects such as metrics, connections, gauge transformations,
curvature, gravity, etc.

\item If considering manifolds in general, one would like to extend the
theory of algebraic topology. Noncommutative geometry may provide a
mechanism for topology change. 

\end{mylist}

In this talk I hope to give an my personal interpretation of the
noncommutative or ``fuzzy'' sphere. The idea is to give the reader
some motivation for looking at my work.  I shall only state the
theorems and the reader is invited to look at my papers for the
details. There are many reasons for studying the sphere.

\begin{mylist}

\item It is simple enough that most of the results can be explicitly
worked out, and being two dimensional and compact, it can also be
pictured.

\item The space of functions on the sphere is a representation of
$su(2)$. As the simplest of classical groups it has been studied in
depth and there are many results which can be used.

\item The sphere has curvature.

\item The results and insight gained from the study of the sphere can
be extended to certain symmetric homogeneous spaces (co-adjoint
orbits) of higher dimensions, and to non compact spaces.

\item It, together with the torus, can be used to study the vector
space, connections and other tools used in differential geometry and
general relativity.

\item It, and the torus are the only 2-dimensional manifolds for which
noncommutative analogues are ``easy''.

\end{mylist}


\section{The Noncommutative Sphere}

In paper \cite{Gratus5} details are given about the noncommutative
analogue of the space of functions on the sphere.  This satisfies (1)
and (2) of the wish list.

\subsection{The commutative sphere: spherical harmonics}

As mentioned in the introduction, we must first consider the algebra,
$C(S^2)$, of complex valued functions from the sphere. Ignoring
problems to do with linear analysis{}\footnote{We avoid the problems
of linear analysis by considering only the algebra of finite sums of
spherical harmonics.}, we can write the functions on the sphere as
sums of {\it Spherical Harmonics} $Y^m_n$, where $n=0,1,2,\ldots$ and
$m=-n,-n+1,\ldots,n$.  In spherical coordinates $(\theta,\phi)$, where
$0\le\theta\le\pi$ and $0\le\phi<2\pi$ we can write the spherical
harmonics in terms of Legendre polynomials:
\begin{align*}
Y^{m}_n(\theta,\phi) 
&=
\left(\frac{(n+m)!(2n+1)}{(n-m)!}\right)^{\scrhalf}
e^{-im\phi} P^{-m}_n(\cos\theta) 
\end{align*}
Alternatively one defines them in terms of a subset of the rotation
matrices.
\begin{align}
Y^{-m}_n(\beta,\alpha) 
&=
\left(\frac{2n+1}{4\pi}\right)^{\scrhalf}
{D^n_{m0}(\alpha,\beta,\gamma)}
\label{Ynm_Dnm}
\end{align}
where $(\alpha,\beta,\gamma)$ are the Euler angles.
(The Spherical Harmonic is, of course, independent of the value of
$\gamma$.)

For our purposes, there is an alternative definition of the spherical
harmonics given in terms of cartesian coordinates on $\Real^3$. 
For a sphere of radius $R$ we define
\begin{align*}
f_{\text{polar}}(\theta,\phi) = f_{\text{cart}}(x,y,z)
\qquad\text{where}\qquad
x=R\cos\phi\,\sin\theta\,,\
y=R\sin\phi\,\sin\theta\,,\
z=R\cos\theta
\end{align*}
It is clear that for any $f_{\text{polar}}(\theta,\phi)$ there will be
many functions $f_{\text{cart}}(x,y,z)$. For example the
functions
\begin{align*}
&f_{\text{cart}}(x,y,z)=x^2+y^2+z^2 
&&\text{and}
&&f_{\text{cart}}(x,y,z)=R^2
\end{align*}
represent the same function on $S^2$. To choose a unique way of
writing the functions on the sphere we require the polynomial to be
{\it symmetric} and {\it formally tracefree}. That is there is a
unique way of writing any{}\footnote{Here any function means any
function which is a nicely convergent sum of spherical harmonics}
function $f:S^2\mapsto\Cmpx$ as
\begin{align}
f(x_1,x_2,x_3) &= f_0 + \sum_{a=1}^3 f_a x_a + 
\sum_{a,b=1}^3 f_{ab} x_a x_b + \ldots +
\sum_{a_1,\ldots a_n=0}^3 f_{a_1 a_2 \cdots a_n} x_{a_1} x_{a_2}\cdots
x_{a_n} + \cdots
\label{fxyz}
\end{align}
where 

\begin{mylist}

\item $f_{a_1 a_2 \cdots a_n}$ is totally symmetric in any two indices.
and

\item $\sum_{b=1}^3 f_{bba_3\cdots a_n}=0$
\end{mylist}

The subspace $\Pexpr^n$ of polynomials of order $n$ 
\begin{align*}
\Pexpr^n &= \text{span}\{f_{a_1 a_2 \cdots a_n}
x_{a_1} x_{a_2}\cdots x_{a_n}\}
\end{align*}
which satisfy the conditions above, has dimension $2n+1$. This space
is spanned by the spherical harmonics
\begin{align*}
Y^{-n}_n\,,\ 
Y^{-n+1}_n\,,\ldots, 
Y^{0}_n\,,\ldots, 
Y^{n-1}_n\,,\ 
Y^{n}_n\,,\ 
\end{align*}
The low order spherical harmonics are given in both polar coordinates
and cartesian coordinates in Table \ref{tab_sph_harm}.
Here $J_+=x+iy$ and $J_-=x-iy$.
\begin{table}[p]
\rotatebox{90}{
\begin{tabular}{|c|l|l|l|}
\hline
Harmonic & 
Polar coordinates (for $\varepsilon=0$) & 
Formally traceless symmetric polynomial ($\forall \varepsilon$) &
Hahn Polynomial ($\forall \varepsilon$)
\\
\hline
$ P^0_0$ & $1$ & $1$ & $1$
\\
\hline
$ P^{1}_1$ &
$Re^{-i\phi}\sin\theta$ &
$J_+$ &
$J_+$ 
\\
$ P^{0}_1$ &
$-\sqrt{2} R\cos\theta$ & 
$-\sqrt{2} z$ &
$-\sqrt{2} z$ 
\\
$ P^{-1}_1$ &
$-Re^{i\phi}\sin\theta$ &
$-J_-$ &
$-J_-$ 
\\
\hline
$P^{2}_2$ &
$R^2e^{-2i\phi}\sin^2\theta$ &
$J_+^2$ &
$J_+^2$ 
\\
$P^{1}_2$ &
$-2R^2e^{-i\phi}\cos\theta\sin\theta$ &
$-2(J_+z + zJ_+)$ &
$-J_+(2z+\varepsilon)$ 
\\
$P^{0}_2$ &
$R^2\sqrt{2/3}(3\cos^2\theta-1)$  &
$\sqrt{2/3}(2z^2 - J_+J_- - J_-J_+)$ &
$\sqrt{2/3}(3z^2-R^2)$
\\
$P^{-1}_2$ &
$2R^2e^{i\phi}\cos\theta\sin\theta$ &
$2(J_-z + zJ_-)$ &
$J_-(2z-\varepsilon)$ 
\\
$P^{-2}_2$ &
$R^2e^{2i\phi}\sin^2\theta$   &
$J_-^2$ &
$J_-^2$ 
\\
\hline
$P^{3}_3$ &
$R^3e^{-3i\phi}\sin^3\theta$ &
$J_+^3$ &
$J_+^3$ 
\\
$P^{2}_3$ &
$-R^3e^{-2i\phi}\sqrt{6}\cos\theta\sin^2\theta$ &
$-\sqrt{6}(J_+^2z+J_+zJ_++zJ_+^2)$ &
$-\sqrt{6}J_+^2(z+\varepsilon)$ 
\\
$P^{1}_3$ &
$R^3e^{-i\phi}\sqrt{3/5}\sin\theta(5\cos^2\theta-1)$ &
\begin{tabular}{@{}l}
$-\sqrt{3/5}(J_+^2J_-+J_+J_-J_++J_-J_+^2)$ \\
\qquad $+4\sqrt{3/5}(J_+z^2+zJ_+z+z^2J_+) $ 
\end{tabular}
&
$\sqrt{3/5}J_+(5z^2+5\varepsilon z + 2\varepsilon^2 - R^2)$ 
\\
$P^{0}_3$ &
$R^3(2/\sqrt{5})(3\cos\theta-5\cos^3\theta)$ &
\begin{tabular}{@{}l}
$(2/\sqrt{5})(z^3
- J_+J_-z - J_+zJ_- - J_-J_+z$\\
\qquad $ - J_-zJ_+ - zJ_+J_- - zJ_-J_+)$ 
\end{tabular} 
&
$(-2/\sqrt{5})z(5z^2-3R^2+\varepsilon^2)$ 
\\
$P^{-1}_3$ &
$-R^3e^{i\phi}\sqrt{3/5}\sin\theta(5\cos^2\theta-1)$ &
\begin{tabular}{@{}l}
$\sqrt{3/5}(J_-^2J_++J_+J_+J_-+J_+J_-^2)$ \\
\qquad $-4\sqrt{3/5}(J_-z^2+zJ_-z+z^2J_-) $ 
\end{tabular}
&
$-\sqrt{3/5}J_-(5z^2-5\varepsilon z + 2\varepsilon^2 - R^2)$
\\
$P^{-2}_3$ &
$-R^3e^{2i\phi}\sqrt{6}\cos\theta\sin^2\theta$ &
$-\sqrt{6}(J_-^2z+J_-zJ_-+zJ_-^2)$ &
$-\sqrt{6}J_-^2(z-\varepsilon)$ 
\\
$P^{-3}_3$ &
$e^{3i\phi}R^3\sin^3\theta$ &
$-J_-^3$ &
$J_-^3$ 
\\
\hline
\end{tabular}
}
\caption{Low order Spherical harmonics written in both polar and
traceless cartesian form (the normalisation may not be conventional)}
\label{tab_sph_harm}
\end{table}

\subsection{The commutative sphere: the inner product}

There is a natural inner product on $C(S^2)$ given by
\begin{align*}
\langle {}{f},{}{g}\rangle_{S^2}
=\frac{1}{4\pi R^2}\int_{S^2} \cnj{{}{f}}{}{g}\sin\theta d\phi
d\theta
\end{align*}
which satisfies
\begin{align*}
\frac1{4\pi R^2}\int Y^m_n(\theta,\phi) \sin\theta d\phi d\theta &=
\left\{ 
\begin{array}{ll}
1 & \qquad m=0 \hbox{ and } n=0 \cr
0 & \qquad \hbox{otherwise}
\end{array} \right.
\nonumber
\end{align*}
therefore we can define the inner product as simply 
\begin{align*}
\innerprod(f,g)_{S^2} = \pi_0(\cnj{f}g)
\end{align*}
where $\pi_0(f)$ means take the $Y^0_0$ coefficient of the $f$, or if
written as (\ref{fxyz}), take the $f_0$ component.


\subsection{The noncommutative algebra $\Pexpr(\varepsilon,R)$}

We construct a two parameter algebra $\Pexpr(\varepsilon,R)$ which may
be thought of as the noncommutative analogue of $C(S^2)$ as follows:
The universal enveloping algebra of $su(2)$ is given by
\begin{align*}
{\cal U}(\varepsilon) &=
\{ \mbox{Free noncommuting algebra of polynomials in $x,y,z$ } \}
\Big/\sim
\end{align*}
where $x,y,z$ obey the commutation relations for $su(2)$
\begin{align}
[x,y] \sim i\varepsilon z,\
[y,z] \sim i\varepsilon x,\
[z,x] \sim i\varepsilon y
\label{Pex_com_rel}
\end{align}
The Casimir is of course given by $x^2+y^2+z^2$. This belongs to the
center of ${\cal U}(\varepsilon)$. We may
therefore quotient ${\cal U}(\varepsilon)$
\begin{align*}
\Pexpr(\varepsilon,R) &=
{\cal U}(\varepsilon) \bigg/ J(R)
\end{align*}
where J(R) is the two sided ideal generated by
\begin{align}
x^2+y^2+z^2 &\sim R^2
\label{Pex_casim}
\end{align}
This gives an infinite dimensional two parameter algebra
$\Pexpr(\varepsilon,R)$ where in the most general case,
$\varepsilon,R\in\Cmpx$.  For different values of $\varepsilon$ and
$R$ we obtain:

\begin{mylist}

\item The commutative algebra of finite sums of harmonics on the
sphere (when $\varepsilon=0$ and $R\in\Real$, $R>0$).  In this case
$R$ is radius of the sphere.

\item The finite matrix representation of $su(2)$.
When $\varepsilon^2(N^2-1)=4R^2$ and $N\in\Intg$, $N\ge1$, then
$M_N(\Cmpx)$ forms a quotient algebra, and $R^2$ is the Casimir
operator.

\item A noncommutative algebra of polynomials which is an infinite
dimensional reducible but non decomposable representation of $su(2)$,
for all values of $\varepsilon\ne 0$.
\end{mylist}

As mentioned, for the case $\varepsilon=0$, the noncommutative algebra
$\Pexpr(\varepsilon=0,R)$ becomes the commutative algebra $C(S^2)$. The
first correction term in the product of two functions is given by the
Poisson bracket on $S^2$, defined with respect to its natural
measure\footnote{It is necessary to map
$C(S^2)\mapsto\Pexpr(R,\varepsilon)$ in order to take the commutator.
This map is well defined.}.
\begin{align*}
\{\bullet,\bullet\}&:\Czz(S^2)\times \Czz(S^2)
\mapsto
\Czz(S^2)
\nonumber
\\
\{ f, g\} &=
\lim_{\varepsilon\to0}
\left(\frac{1}{i\varepsilon} 
[f,g] \right)
\\
\{ f, g\} &= 
\frac1{R\sin\theta}\left(
\frac{\partial  f}{\partial\phi} \, 
\frac{\partial  g}{\partial\theta} 
-
\frac{\partial  f}{\partial\theta} \, 
\frac{\partial  g}{\partial\phi} 
\right)
\end{align*}
It is still unknown (at least to me) the exact nature of the higher
order correction terms, and whether or not they form a $\star$-product
in the sense of Flato et al. I am looking into this problem.

\vskip 1em

We can still represent the functions using (\ref{fxyz}). Thus we still
have a sesquilinear form on $\Pexpr(\varepsilon,R)$ given by 
\begin{align*}
\innerprod(f,g)=\pi_0(f^\dagger g)
\end{align*}
where $\pi_0(f)$ means take the $f_0$ component of $f$, and 
$f^\dagger$ is the hermitian conjugate of $f$ defined by:
\begin{align*}
\dagger:\Pexpr\mapsto\Pexpr,\
(ab)^\dagger = b^\dagger a^\dagger,\ 
x^\dagger=x,\ 
y^\dagger=y,\ 
z^\dagger=z,\ 
\lambda^\dagger=\cnj\lambda \qquad\hbox{for $\lambda\in\Cmpx$}
\end{align*}
However this sesquilinear form is no longer positive definite since
$\pi_0(f,f)$ may be positive, negative or zero. (See (\ref{Norm_Pmn})
below.) 

\subsection{$P^m_n$ an orthogonal basis for $\Pexpr$}

There is a basis\footnote{In \cite{Gratus5} $P^m_n$ are normalised by
using a constant $\alpha_n$. Also we use $\kappa$ in place of
$\varepsilon$.} of $\Pexpr(\varepsilon,R)$ given by
\begin{align}
\{P^m_n\,|\ n,m\in\Intg,n\ge0,|m|\le n\}
\end{align}
where (for $\varepsilon\ne0$)
\begin{align}
P^m_n &=
\varepsilon^{m-n}
\left(\frac{(n+m)!}{(2n)!\,(n-m)!}\right)^{\scrhalf} 
\left(\Ad{J_-}\right)^{n-m}(J_+{}^n)
\label{Pmn_def_Pmn}
\end{align}
where $J_\pm=x\pm iy$ and $\Ad{J_-}f=[J_-,f]$.  This basis is
orthogonal with respect to the sesquilinear form
\begin{align}
\innerprod({P^{m_1}_{n_1}},{P^{m_2}_{n_2}})
&=
\delta_{m_1 m_2}
\delta_{n_1 n_2}
\norm{P^{m_1}_{n_1}}^2
\end{align}
where the norm of $P^m_n$ is independent of $m$ and is
given by
\begin{align}
\norm{P^m_n}^2 &=
\frac{(n!)^2}{(2n+1)!} 
\prod_{r=1}^n(4R^2+\varepsilon^2(1-r^2))
\label{Norm_Pmn}
\end{align}
We can see that $\norm{P^m_n}^2$ may be positive, negative, or zero
depending on the values od $R,\varepsilon$ and $n$.

Each $P^m_n$ can be written as a homogeneous formally tracefree
symmetric polynomial in $(x,y,z)$ of order $n$. As such they are
independent of the value of $\varepsilon$ and $R$. Therefore they are
exactly the same as the spherical harmonics when written in this form,
given in table \ref{tab_sph_harm}.

Each $P^m_n$ is an eigenvector of the operators $\Ad{z}$
and $\Delta=\Ad{z}^2+\Ad{J_+}\Ad{J_-}+\Ad{J_-}\Ad{J_+}$.
\begin{align}
\Ad{z} P^m_n &= \varepsilon m P^m_n 
\label{Pmn_ez_Pmn} 
\\
\Delta P^m_n &= \varepsilon^2 n(n+1) P^m_n 
\label{Pmn_Del_Pmn} 
\end{align}
The ladder operators $\Ad{J_+},\Ad{J_-}$ increase and decrease $m$ so
that $\Pexpr^n$ can be viewed as a $2n+1$ dimensional adjoint
representation of $su(2)$.
\begin{align}
\Ad{J_+} P^m_n &= 
\varepsilon (n-m)^{\scrhalf} (n+m+1)^{\scrhalf} P^{m+1}_n 
\label{Pmn_ep_Pmn}
\\
\Ad{J_-} P^m_n &= 
\varepsilon (n+m)^{\scrhalf} (n-m+1)^{\scrhalf} P^{m-1}_n 
\label{Pmn_em_Pmn}
\end{align}
The effect of taking the hermitian conjugate is given by
\begin{align}
(P^m_n)^\dagger &= (-1)^m P^{-m}_n 
\label{Pmn_Pmn_dag}
\end{align}
The formula for the combination of 2 basis elements may be obtained by
substituting $r_1=r_2=r=0$ and
$\varepsilon^2(k-\tfrac12)^2=R^2+\tfrac14$ into (\ref{RM_res}) below.


\subsection{Alternative ways of writing $\Pexpr$ and $P^m_n$}

As well as writing the elements of $\Pexpr$ as formally traceless
symmetric polynomials, there are at least three other ways of writing
them. As stated for certain $R$ and $\varepsilon$ they may be
quotiented to form a matrix algebra.  We may also write the elements
of $\Pexpr$ in terms of Hahn polynomials and in terms of stereographic
projection (see section \ref{ch_stereo}). These alternative
representations are useful for some of the proofs as well as leading
to a greater understanding.

\subsubsection{Finite dimensional representations}

If $R^2=\varepsilon^2(k(k+1))$ then we may quotient out all basis
elements $P^m_n$, with $n\ge 2k+1$. This leaves just the $(2k+1)^2$
elements $\{P^m_n | n\le 2k\}$. This algebra is exactly equivalent to
the $(2k+1)\times(2k+1)$ matrix representation of $su(2)$.
\begin{align*}
J_0\rvec(k,j) &= \varepsilon j \rvec(k,j)  
\nonumber\\
J_+\rvec(k,j) &= \varepsilon(k-j)^\scrhalf(k+j+1)^\scrhalf\rvec(k,j+1)
\nonumber\\
J_-\rvec(k,j) &= \varepsilon(k+j)^\scrhalf(k-j+1)^\scrhalf\rvec(k,j-1)
\end{align*}
This representation may be written in terms of Wigner's operators, see
(\ref{RM_Wig_op}) below. The sesquilinear form is now positive
definite\footnote{From (\ref{Norm_Pmn}) one see that
$\norm{P^m_n}^2>0$ for $n<2k+1$ and that $\norm{P^m_n}^2=0$
otherwise}, and given by the trace.
\begin{align*}
\innerprod(f,g) &=
\frac1{2k+1}\sum_{j=-k}^k \lvec(k,j)f^\dagger g\rvec(k,j)
\end{align*}

\subsubsection{Writing $P^m_n$ in terms of Hahn Polynomials}

Given $f\in\Pexpr$ then we can use the commutation relations
(\ref{Pex_com_rel}) to push the $J_+$ and $J_-$ to the left of each
term. If a $J_+$ and $J_-$ appear in one term we can use the Casimir
(\ref{Pex_casim}) identity to remove both. Thus the resulting terms
must either have only $J_+$'s or only $J_-$'s or neither.  If we
collect all the terms with the same number of $J_+$ or $J_-$ as their
factors then $f$ may  be written as a sum of terms of the form
\begin{align}
\{(J_+)^m p(z)\,,\, p(z)\,,\, (J_-)^{-m}p(z)\}
\nonumber
\end{align}
where $p(z)$ is a polynomial in $z$. 
Since $P^m_n$ is a eigenvector of $\Ad{z}$ we have
\begin{align*}
P^m_n = \left\{
\begin{array}{ll}
(J_+)^m
(-\varepsilon)^{n-m} \ppmatrix{2n \cr n-m}^{-\scrhalf}
h^{(m,m)}_{n-m}(\nfrac{z}{\varepsilon}+\tfrac{N-1}2,N-m)
& \hbox{\rm\ for }  m\ge 0 
\\
(J_-)^{-m} 
(-1)^m(\varepsilon)^{n-m} \ppmatrix{2n \cr n-m}^{-\scrhalf}
h^{(m,m)}_{n-m}(\nfrac{z}{\varepsilon}-m+\tfrac{N-1}2,N-m)
& \hbox{\rm\ for }  m<0 
\end{array}
\right.
\end{align*}
where $N^2+1=4R^2\varepsilon^{-2}$.  Here $h^{(\alpha,\beta)}_n(x,N)$
is a Hahn Polynomial\footnote{The Hahn polynomials are defined in a
similar manner to Legendre polynomials, except one replaces the
integral by a finite sum in the definition of the orthogonality
property. In the limit $N\to\infty$ they tend to the Legendre
polynomials}, following the notation of \cite[chapter 2]{Nik}.
Examples of $P^m_n$ for $n\le 3$ written in this form are given in
table \ref{tab_sph_harm}.


\section{The ``Wick rotation'' of $\Pexpr$ to the Hyperbolic case}

The space of scalar functions on the sphere is a representation of
$su(2)$. Now we consider extending this for a representation of
$sl(2,\Cmpx)$. This algebra contains both $su(2)$ and $su(1,1)$ as
subalgebras. We can extend the results for sphere to the
symmetric spaces of $su(1,1)$ viz the positive hyperboloid
(hyperbolic disc) and
negative hyperboloid, (two dimensional De-Sitter space).

Let us consider a new version of the algebra of noncommutative
scalar functions given by: 
\begin{align}
\Pexpr(\varepsilon,R,\alpha) &= \left\{
\mbox{Polynomials in $J_+,J_-,J_0$}\right\} \bigg/ =
\end{align}
where
\begin{align}
[J_0,J_+] &= \varepsilon J_+ &
[J_0,J_-] &= -\varepsilon J_- &
[J_+,J_-] &= 2\varepsilon\alpha^2 J_0 &
J_0^2 + \frac1{2\alpha^2}(J_+J_- + J_-J_+) &= R^2
\label{sl_com_rel}
\end{align}
In the most general case $\varepsilon,R,\alpha\in\Cmpx$ are all
independent and $\alpha\ne0$.  The hermitian conjugate is given by
$J_0^\dagger=J_0$, $J_+^\dagger=J_-$ and $J_-^\dagger=J_+$, and the
sesquilinear form is defined as before. By extending the results
of \cite{Gratus5} we have very similar results for the new basis
elements $P^m_n$ except with $\alpha$'s dotted about.
A significant result is that for $\alpha^2=-1$ and
$R^2<-\varepsilon^2$, the sesquilinear form of
$\Pexpr(\varepsilon,R,\alpha)$ positive definite, making
$\Pexpr(\varepsilon,R,\alpha)$ into a Hilbert space. However the
operations on this spaces of multiplying by the elements
$J_0,J_+,J_-$ are not continuous. There must be some way of recovering
the representation of $su(1,1)$. This is being investigated. 

By using the Jordan-Schwinger representation of $sl(2,\Cmpx)$ we should
be able to generate $\baseterm(n,r,m)$, to produce a theory of spinor
and vector fields on the positive and negative hyperboloids.


\subsection{Stereographic projection for $\Pexpr$}
\label{ch_stereo}

The algebra $\Pexpr$ is equivalent to the algebra given
by\footnote{here $z$ and $x$ are not the same as the $z$ in section 2,
$z$ is replaced by $J_0$ and $x$ by $\frac12(J_++J_-)$}
\begin{align}
\Pexpr_z &= \left\{ 
\mbox{Polynomials in } z,\cnj z, \frac1{\Rp^2+\alpha^2 z\cnj z} 
\right\} \bigg/ \sim
\end{align}
where $z$ and $\cnj z$ are considered as independent and conjugate
$z^\dagger=\cnj z$, and $\Rp^2=R^2+\tfrac14\varepsilon^2$ as before.
The quotient is given by
\begin{align}
z\cnj z - \cnj z z &\sim
\frac{-\varepsilon}{8\Rp^3\alpha^2}
(4\Rp^2 + \alpha^2 z\cnj z)
(4\Rp^2 + \alpha^2 \cnj z z)
\end{align}
For the case $R=1$ and $\alpha^2=-1$, this reduces to the formula
given in \cite{Klimek_Lesn1}. This is equivalent to
\begin{align}
\cnj z z &= \rho(z\cnj z)
\end{align}
where $\rho$ is the m\"obius transform
\begin{align}
\rho(x) &= \frac{
(1+{\varepsilon}/{2\Rp})x + 2\varepsilon\Rp/\alpha^2}
{(-\varepsilon\alpha^2/8\Rp)x + (1-{\varepsilon}/{2\Rp})}  
\end{align}
This m\"obius transform when written as a function of both
$\varepsilon$ and $x$ satisfies $\rho^n(\varepsilon,x) =
\rho(n\varepsilon,x)$ where
$\rho^{n+1}(\varepsilon,x)=\rho(\varepsilon,\rho^n(\varepsilon,x))$.

We may consider this the analogue of stereographic projection. The
relationship between the two algebras is given by (where $x=z\cnj z$)
\begin{align}
J_0 &= \Rp\frac{4\Rp^2 - \alpha^2 x}{4\Rp^2 + \alpha^2 x} -
\frac{\varepsilon}2
&
J_+ &=
i\cnj z
\frac{4\Rp^2\alpha^2}{4\Rp^2 + \alpha^2 x}
&
J_- &=
-i
\frac{4\Rp^2\alpha^2}{4\Rp^2 + \alpha^2 x} z
\end{align}

We note that if $f\in\Pexpr^m$ then the denominator of $f$ is
$(\Rp^2+\alpha^2 x)^n$, and the numerator is a polynomial of $z$,
$\cnj z$ and $x$. Thus the numerator of $P^m_n$ will be a polynomial
in $z$ and $x$ or $\cnj z$ and $x$ depending on the sign of $m$. This
polynomial is related to the Hahn polynomials and its degree will be
less than $n$.

There also exists projections for mapping the noncommutative negative
hyperboloid to a noncommutative cylinder.


\subsection{The commutative limit}

When $\varepsilon=0$ the algebra $\Pexpr$ reduces to the commutative
algebra of functions on either the sphere or the hyperboloid depending
on the sign of $\alpha^2$.  For the sphere and the positive
hyperboloid the algebra $\Pexpr_z$ become the stereographic
projection, whilst for the negative hyperboloid we have to consider
the analytic continuation of the stereographic projection. In the case
of the cones the stereographic projection is degenerate. See table for
the four symmetric spaces of $su(2)$ or $su(1,1)$.

\begin{table}[t]
\begin{tabular}{|@{}l@{}|@{}c@{}|@{}c@{}|@{}c@{}|@{}c@{}|}
\hline
& Sphere & Positive hyperboloid & Negative hyperboloid & Two cones 
\\
&& Hyperbolic disc & 2 dim De-Sitter space &
\\
&
\epsfysize=3 cm
\epsffile{circle.eps}
&
\epsfysize=3 cm
\epsffile{poshyp.eps}
&
\epsfysize=3 cm
\epsffile{neghyp.eps}
&
\epsfysize=3 cm
\epsffile{twocone.eps}
\\
\hline
&
$\alpha^2=1$ &
$\alpha^2=-1$ &
$\alpha^2=-1$ &
$\alpha^2=-1$ 
\\
&
$R^2>0$ &
$R^2>0$ &
$R^2<0$ &
$R^2=0$ 
\\
\begin{tabular}[t]{l}
coordinate \\ system 
\end{tabular}
&
$\displaystyle{
\begin{aligned}[t]
J_0&=R\cos\theta \\
J_+&=ie^{-i\phi}R\sin\theta \\
J_-&=-ie^{-i\phi}R\sin\theta \\
0 &\le \theta \le\pi \\
0 &\le \phi < 2\pi 
\end{aligned} }$
&
$\displaystyle{
\begin{aligned}[t]
J_0&=R\cosh\eta \\
J_+&=ie^{-i\phi}R\sinh\eta \\
J_-&=-ie^{-i\phi}R\sinh\eta \\
0 &\le \eta < \infty \\
0 &\le \phi < 2\pi 
\end{aligned} }$
&
$\displaystyle{
\begin{aligned}[t]
J_0&=i{R}\sinh\eta \\
J_+&=-e^{-i\phi}{R}\cosh\eta \\
J_-&=e^{-i\phi}{R}\cosh\eta \\
-\infty &< \theta < \infty \\
0 &\le \phi < 2\pi 
\end{aligned} }$
&
$\displaystyle{
\begin{aligned}[t]
J_0&=\eta \\
J_+&=ie^{-i\phi}|\eta| \\
J_-&=-ie^{-i\phi}|\eta| \\
-\infty &< \eta < \infty \\
0 &\le \phi < 2\pi 
\end{aligned} }$
\\
\hline
\begin{tabular}[t]{l}
Stereographic \\
coordinates \\
$z=-ire^{i\phi}$, \\
$z=ire^{-i\phi}$
\end{tabular}
&
\begin{tabular}[t]{l}
$r=2R\tan(\theta/2)$ \\
$z\in\Cmpx$  
\end{tabular} 
&
\begin{tabular}[t]{l}
$r=2R\tanh(\theta/2)$ \\
$z\in\Cmpx$, $|z|<1$   
\end{tabular} 
&
\begin{tabular}[t]{l}
$r=2iR\tanh(\eta/2-i\pi/2)$ \\
$z$ is not a complex \\
variable 
\end{tabular} 
&
degenerate
\\
\hline
\end{tabular}
\caption{The two dimensional symmetric spaces of $sl(2,\Cmpx)$}
\label{tbl_comm_spaces}
\end{table}


\section{Spinors and Vectors on the Noncommutative Sphere}

In \cite{Gratus6} we consider how one can extend the algebra
$\Pexpr(\varepsilon,R)$ to include vector and spinor fields. The
result is a new algebra $(\PsiSpace,\rho)$, where $\PsiSpace$ is a set
of polynomials, and $\rho$ is a noncommutative and nonassociative
product on $\PsiSpace$. There is a natural basis of $\PsiSpace$ given
by $\{\baseterm(n,r,m)\}$ where $n=0,\pm\frac12,\pm1,\pm\frac32\ldots$
and $m,r=-n,-n+1,\ldots,n$. It is useful to consider the subspace
\begin{align*}
\PsirSpace{r}=\text{span}\{\baseterm(n,r,m)\ |\ \forall n,r\}
\end{align*}
The subspace $\PsirSpace{0}$ together with the product $\rho$ may be
identified with $\Pexpr(\varepsilon,R)$, the space of scalar functions
on the sphere. Thus $\baseterm(n,0,m)=P^m_n$. This is analogous to
(\ref{Ynm_Dnm}).  All other $\PsirSpace{r}$ are modules over
$\PsirSpace{0}$ in an analogous way to the way spinors and vector
fields are modules over scalar fields.  The spaces
$\PsirSpace{1},\PsirSpace{-1},\PsirSpace{1/2}$ will be identified with
vector, covector, and spinor fields respectively.

In the limit, $\varepsilon=0$, the basis element $\baseterm(n,r,m)$
becomes the rotation matrix element $D^n_{mr}(\alpha,\beta\gamma)$.
As a result we may view $(\PsiSpace,\rho)$ as the noncommutative and
nonassociative analogue of the algebra of functions on the Lie group
$SU(2)$.


\subsection{The mechanics of setting up $\PsiSpace$}

The first two chapters of \cite{Gratus6} concern themselves with the
setting of $(\PsiSpace,\rho)$.  For an outline, consider the
Jordan-Schwinger representation of $su(2)$ given by
\begin{align*}
J_0 &= \tfrac12(a_+a_--b_+b_-) &
J_+ &= a_+b_- &
J_- &= a_-b_+
\end{align*}
where $[a_-,a_+]=\varepsilon$ and $[b_-,b_+]=\varepsilon$ are the
generators of $\Weil$ the product of two Heisenberg-Weil algebras. The
Casimir of this representation of $su(2)$ is given by
\begin{align*}
J_0^2 + \tfrac12 J_+J_- + \tfrac12 J_-J_+ &= K_0^2 - 
\tfrac14 \varepsilon^2 
\end{align*}
where 
\begin{align*}
K_0 &= \tfrac12(a_+a_-+b_+b_-+\varepsilon)
\end{align*}
We would like to take the square root of this equation and consider
quotienting $\Weil$ by the ideal generated by
$K_0\sim\Rp=(R^2+\tfrac14\varepsilon^2)^\scrhalf$. However, although
$K_0$ commutes with any polynomial in $\{J_0,J_+,J_-\}$, it does not
commute with all the elements in $\Weil$. Therefore the left ideal
generated by $K_0\sim\Rp$ is not a two sided ideal and the
corresponding quotient product is not associative.

Many results about the definition and basis elements
$\baseterm(n,r,m)$ are similar to the corresponding case for the
scalar functions $P^m_n$, but with caution reflecting the
nonassociative nature of $(\PsiSpace,\rho)$.  They are not given here
and the reader is invited to look at the article.

However, I will quote the product formula for the basis elements,
since it is not given in the above paper, but will hopefully appear in
\cite{Gratus7}. The combination of two basis elements is given by
\begin{align}
\rho(\baseterm(n_1,r_1,m_1) 
\baseterm(n_2,r_2,m_2) )
&=
\sum_{n=n_{\min}}^{n=n_1+n_2}
\CG(n_1,n_2,n,m_1,m_2,m_1+m_2)
\RMpsi(n_1,n_2,n,r_1,r_2,r_1+r_2)
\baseterm(n,r_1+r_2,m_1+m_2) 
\label{RM_def}
\end{align}
where $n_{\min}=\max(|n_1-n_2|,|r_1+r_2|,|m_1+m_2|)$,
$\CG(n_1,n_2,n,m_1,m_2,m_1+m_2)$ is the Clebsh-Gordon coefficient, and
the reduced matrix element $\RMpsi(n_1,n_2,n,r_1,r_2,r_1+r_2)$ is
given by\footnote{Here
$\norm{\baseterm(n,r,\bullet)}_k$ is given by substituting
$R=\varepsilon(k(k+1))^\scrhalf$ into the formula for
$\norm{\baseterm(n,r,\bullet)}^2$ and taking the positive square
root.}
\begin{align}
\RMpsi(n_1,n_2,n,r_1,r_2,r_1+r_2)
&=
(-1)^{2k+n_1+n_2+r_1+r_2}
\frac{
\norm{\baseterm(n_1,r_1,\bullet)}_{k+r_2}
\norm{\baseterm(n_2,r_2,\bullet)}_{k}
}{\norm{\baseterm(n,r_1+r_2,\bullet)}_k}
(2k+2r_2+1)^\scrhalf \times
\notag\\
&\qquad\qquad\qquad\qquad\qquad\qquad\qquad
(2n_1+1)^\scrhalf (2n_2+1)^\scrhalf 
\WsixJ(k+r_1+r_2,n_1,k+r_2,n_2,k,n) 
\label{RM_res}
\end{align}
where the symbol in the curly brackets is Wigner's 6-$j$ coefficient.

This is proved by considering $\baseterm(n,r,m)$ as one of Wigner's
operators
given in \cite[eqn (3.340)]{BL1} as
\begin{align}
\baseterm(n,r,m) \rvec(k,j)
&=
(-1)^{n-r} 
\norm{\baseterm(n,r,\bullet)}_k
\frac{(2n+1)^\scrhalf (2k+1)^\scrhalf}{(2k+2r+1)^{\scrhalf}}
\Wigner(2n,n+r,n+m)
\rvec(k,j)
\label{RM_Wig_op}
\end{align}
and then using the product law given by \cite[eqn (3.350)]{BL1}.  In
\cite{BL1} the Wigner's operators are indeed referred to as discrete
rotation matrices.


\subsection{Physical interpretation}

As mentioned in the introduction, the basis elements can be viewed as
the nonassociative and noncommutative analogue of the rotation
matrices.  However they may also be interpreted in terms of vectors
and spinors on the sphere. This interpretation is not as clear cut as
the case for the scalars and I hope it is not the final word on the
matter.


\subsubsection{The exterior derivative, 1-forms, and vector fields}
\label{ch_covec}

In standard geometry there are many equivalent definitions of a vector
field. However not all of them can be extended to the noncommutative
case at the same time. If ones take the definition of a vector field
given that it must satisfy the Leibniz formula, and extends this
definition to noncommutative geometry then the vector fields are given
by $X=\Ad{f}$ for some $f\in\Pexpr=\PsirSpace{0}$.
This can then be used to give a
definition of the exterior derivative and the exterior algebra
\cite{Madore_bk}.  The problem with this definition is
that this space does not form a module over the space of functions.

Here we give another definition of vectors which do form a module over
the space of functions. We start by defining the covectors by use of
the analogue of the following definition of the exterior derivative of
scalar fields:
\begin{align}
df &= \sum_i \frac{\partial f}{\partial x^i} dx^i
\label{Vec_df_partial}
\end{align}
If we consider the noncommutative sphere as a three dimensional
manifold with normalised basis coordinates
$x^m=(2\Rp+\varepsilon)^{-\scrhalf}\baseterm(1,0,m)$ for
$m\in\{-1,0,1\}$ we can define the basis 1-forms as
\begin{align*}
dx^m=(2\Rp+\varepsilon)^{-\scrhalf}\baseterm(1,-1,m)
\qquad\text{for } m\in\{-1,0,1\}
\end{align*}
This enables us to define the \defn{exterior derivative} on the space
of functions as
\begin{align}
&d:\PsirSpace{0}\mapsto\PsirSpace{-1}
&
df &= \sum_{m=-1}^1 (-1)^{m+1} dx^{-m} \Ad{x^m} f
\label{Vec_def_df}
\end{align}
which is analogous to (\ref{Vec_df_partial}). This is not a derivative
but does satisfy the Leibniz rule in the limit as $\varepsilon\to0$
\begin{align}
d(fg)=d(f)g + f\, d(g) + O(\varepsilon)
\label{Vec_d_limit}
\end{align}

We can now define the {\it vectors fields} as the set $\PsirSpace{1}$.
The action of a vector on a scalar is given by
\begin{align*}
X (f) &= \rho((df)\,X)
\end{align*}
By the same reasoning as (\ref{Vec_d_limit}), this is also only a
derivation in the limit $\varepsilon\to0$.
\begin{align*}
X(fg) &= X(f) g + f X(g) + O(\varepsilon)
\end{align*}
One must therefore decide which of the properties of vector fields one
wishes to extend to noncommutative geometry, since the property of
being a derivative, and the property of being a module over the space
of functions are incompatible. 

We also observe that
\begin{align*}
\sum_{m=-1}^1 x^m X_m = \sum_{m=-1}^1 X_m x^m =  0
\end{align*}
identically for all $\varepsilon$, where $X_m=(dx^m)^\dagger$.
This is equivalent to requiring that $x^m (\partial/\partial x^m)=0$,
i.e. $X^m$ are vector fields on a sphere.

We can now define a \defn{metric}
\begin{align*}
&g:\PsirSpace{1}\times\PsirSpace{1}\mapsto\PsirSpace{0}
&& g(X,Y) = \rho(X^\dagger Y)
\end{align*}
This can be extended for all spaces $\PsirSpace{r}$. 

\subsubsection{Problems with this interpretation}

There are many problems with this interpretation of our system. Here
is a list of cases where this noncommutative geometry is different from
the commutative geometry even in the limit $\varepsilon=0$:

(1) the space $\PsirSpace{-1}$ is the image of $\PsirSpace{0}$ under
$d$. This means that all 1-forms are closed. 

(2) The definition (\ref{Vec_def_df}) can be extended as a map 
$d:\PsirSpace{r}\mapsto\PsirSpace{r-1}$ for all $r$. However this map
does not satisfy $d^2=0$ even in the limit. As a result no exterior
calculus is defined here.

(3) The set of vectors $X_i$ are not dual, in the usual sense, to the
set of 1-forms $d x^j$, or alternatively the vectors $X_i$ are not
orthogonal with respect to the metric $g$. This is because for $i\ne j$
\begin{align*}
g(X_{i},X_{j}) = \rho(dx^{i} X_{j}) = X_{j}(x^{i}) \ne 0 
\end{align*}
for all $\varepsilon$ even when $\varepsilon=0$. However we do have
$\pi_0(g(X_{i},X_{j}))=\delta_{ij}$ 

Some of these problems, together with the non derivative nature of
$d$, may be solved by redefining the space of covectors. For example
it may be similar to (\ref{Spin_def_mat}) below.

\subsubsection{Spinor fields}
\label{ch_spin}

It is natural now to define the set $\PsirSpace{1/2}$ as the space of
spinor fields and its dual $\PsirSpace{-1/2}$. This means that a
vector is the product of two spinors.  Using the usual definition of
rotation by $2\pi$, we see that rotation by $2\pi$ does not change the
sign of $\baseterm(n,r,m)$ if $r$ is an integer, i.e. for scalars,
vectors and other ``Bosons''.  Whilst $\baseterm(n,r,m)$ changes sign
under rotation of $2\pi$ for $r$ a half integer, i.e. for
spinors fields and other ``Fermions''.

An alternative way to define spinor fields 
is as the set
\begin{align}
\spinset &= 
\left\{\binom{a_+}{b_+}f_1 + \binom{a_-}{b_-}f_2 \ \bigg|\
f_1,f_2\in\PsirSpace{0}\right\}
\label{Spin_def_mat}
\end{align}
This is decomposed into $\spinset=\spinset_+\oplus\spinset_-$
corresponding to the eigenspaces of $\Ad{K_0}$, which is now regarded
as the chirality operator.  This interpretation is now equivalent to
the one proposed by Grosse et al \cite{Grosse5,Grosse3,Grosse6}, who
go on to define and solve the Dirac equation.  We note that $\spinset$
is not equivalent to simply two copies of
$\PsirSpace{-1/2}\oplus\PsirSpace{1/2}$, but a proper subset. This is
because the element $\displaystyle{\binom{a_+}{0}\not\in\spinset}$.
In the commutative limit $\varepsilon=0$ we can obtain the standard
spinors on a sphere.


\section{Outlook: A Path Towards Quantum Gravity}

\begin{figure}[t]
\setlength{\unitlength}{1 cm}
\begin{picture}(17.5,13)(0,2)


\put(8.5,14){\framebox(4,1)%
{\bititem\label{dg_fn_S2} 
Functions on $S^2$ }}

\put(4.5,11){\framebox(3.5,2){\parbox{3cm}%
{\bititem\label{dg_Vec_S2} 
Vectors and Spinors on $S^2$. The supersphere}}}

\put(9,11){\framebox(3.5,2){\parbox{3cm}%
{\bititem\label{dg_fn_Kah} 
Functions on compact symmetric manifolds}}}

\put(13.5,11){\framebox(3.5,2){\parbox{3cm}%
{\bititem\label{dg_fn_Lob}
 Functions on +ve and -ve hyperboloid}}}

\put(0,11){\framebox(3.5,2){\parbox{3cm}%
{\bititem\label{dg_connect}
Algebraic \\ connections}}}

\put(1.5,9){\framebox(4,1)
{\bititem\label{dg_connect_S2} 
Connections on $S^2$}}

\put(10.5,9){\framebox(4,1.5){\parbox{3.5cm}%
{\bititem\label{dg_fn_M4} 
Functions on 4 dim De-Sitter space}}}

\put(10.5,6.5){\framebox(4,1.5){\parbox{3cm}%
{\bititem\label{dg_Vec_M4} 
Vectors and  \\ Spinors on De Sitter}}}

\put(3,5.5){\framebox(5,1.5){\parbox{4.5 cm}%
{\bititem\label{dg_connect_M4}
Connections on De Sitter Quantum Gravity}}}

\put(9.5,3.5){\framebox(3.5,2){\parbox{3cm}%
{\bititem\label{dg_diff_stuc}
Existence of differential structure on non Lie algebras}}} 

\put(14,3.5){\framebox(3.5,2){\parbox{3cm}%
{\bititem\label{dg_fn_genus}
Functions on Manifolds of higher genus}}}

\put(3,2){\framebox(9,1){
{\bititem\label{dg_QG_top}
Quantum Gravity with topology change}}}

\put(10.5,14){\vector(-4,-1){3.8}}
\put(10.5,14){\vector(0,-1){0.9}}
\put(10.5,14){\vector(4,-1){3.8}}
\put(2,11){\vector(1,-1){0.9}}
\put(6,11){\vector(-1,-1){0.9}}
\put(11,11){\vector(2,-1){0.9}}
\put(16,11){\vector(0,-1){5.4}}
\put(15,11){\vector(-2,-1){0.9}}
\put(12.5,9){\vector(0,-1){0.9}}
\put(7,11){\vector(1,-1){3.2}}
\put(3.5,9){\vector(1,-1){1.8}}
\put(10.5,7){\vector(-3,-1){2.3}}
\put(5.5,5.5){\vector(1,-1){2.3}}
\put(11,3.5){\vector(0,-1){0.3}}
\put(15.5,3.5){\vector(-4,-1){3.5}}

\end{picture}

\caption{My path towards Quantum gravity}
\label{dg_diagram}

\end{figure}
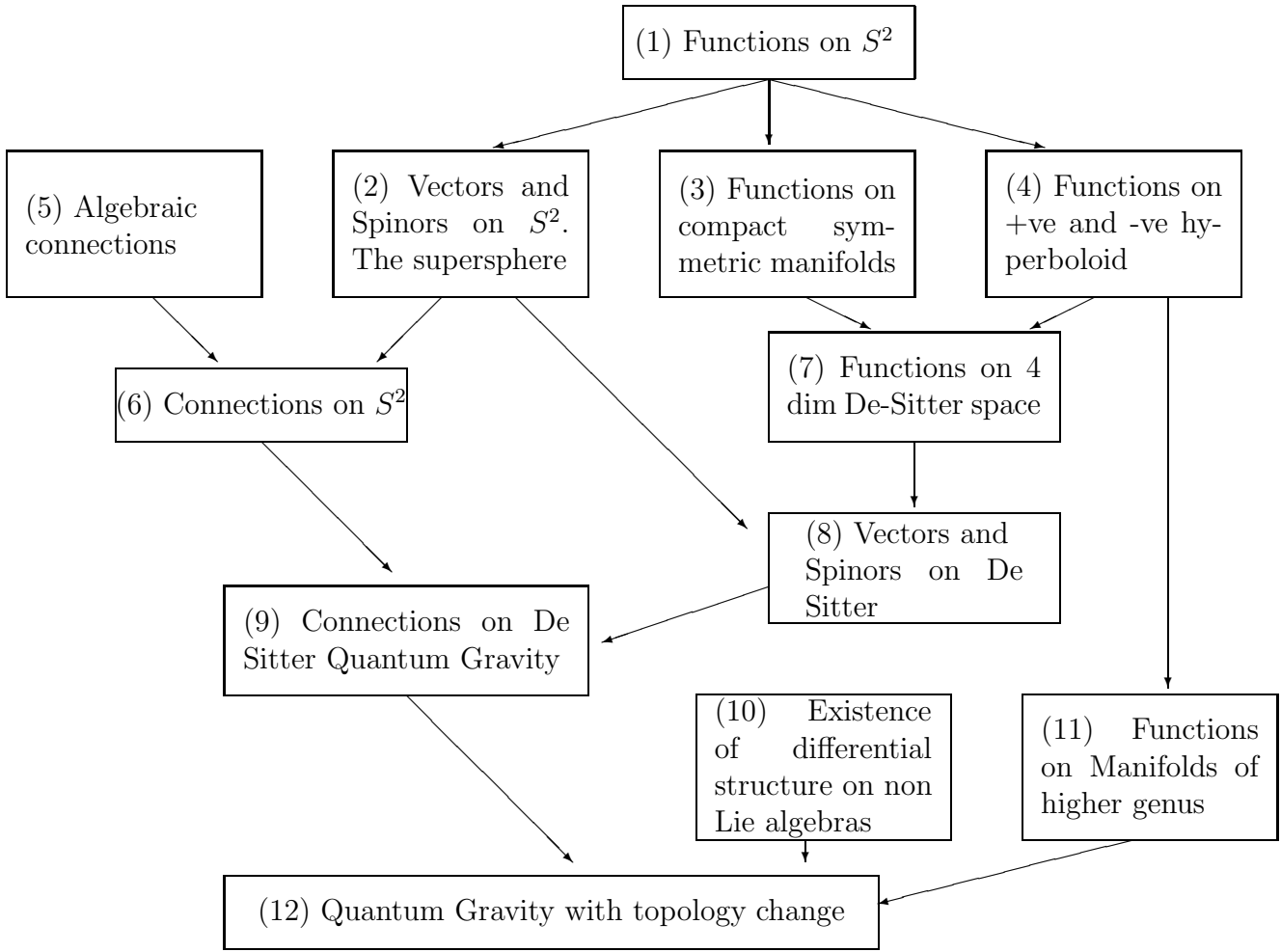


Figure \ref{dg_diagram} shows how the articles described fit into
my personal quest for quantum gravity. 

\begin{list}{}{}

\item[(\ref{dg_fn_S2})]
Paper \cite{Gratus5} as described in section 2 has three natural
extensions $\{(\ref{dg_Vec_S2}),(\ref{dg_fn_Kah}),(\ref{dg_fn_Lob})\}$

\item[(\ref{dg_Vec_S2})]
Paper \cite{Gratus6} as described in section 3. One can also consider
the noncommutative supersphere, which is a representation of the
superalgebra $osp(2,1)$. Thus one can study possible analogies of 
the Dirac operator. 

\item[(\ref{dg_fn_Kah})] 
It was possible to give a noncommutative version of the sphere because
it was a co-adjoint orbit. By using other groups one could extend
these results to higher dimensional compact co-adjoint orbits.  One
could consider the algebra generated by the quotienting the enveloping
algebra of a Lie algebra by the set of equations equivalent to the
Casimirs. As a long term goal one must consider noncommutative
analogous of non co-adjoint orbits.

\item[(\ref{dg_fn_Lob})]
Paper \cite{Gratus7} (hopefully) as described in section 4.

\item[(\ref{dg_connect})] 
The work by Madore \cite{Madore_bk} and Coquereaux \cite{Coquereaux}
examine an algebraic definition connection which can be used in
noncommutative geometry. It requires simply that the vector space of
vector fields is a module over the space of functions. This is the
case in (\ref{dg_Vec_S2}).

\item[(\ref{dg_connect_S2})]
As stated the work in (\ref{dg_connect}) is applicable to
(\ref{dg_Vec_S2}).

\item[(\ref{dg_fn_M4})]
By combining (\ref{dg_fn_Kah}) and (\ref{dg_fn_Lob}) we can generate
an algebra for the 4 dimensional De-Sitter space. This
would give this spacetime a natural cellular structure and enable one
to study scalar fields on this space.

\item[(\ref{dg_Vec_M4})]
The Jordan-Schwinger algebra can be generated for any Lie algebra. It
may be possible therefore to extend the results of (\ref{dg_Vec_S2})
to generate an algebra for spinors and vector fields on the 4
dimensional De-Sitter space.  Since this will enable one to look at
spinor and vector fields on a space with a natural cellular structure,
it might give an alternative method of renormalisation.

\item[(\ref{dg_connect_M4})]
It may be possible to put a (non metric) connection on the
noncommutative 4-dimensional De-Sitter universe and thus have a
version of quantum gravity.

\item[(\ref{dg_diff_stuc})]
Work by myself \cite{Gratus4} and Madore et al \cite{Madore_Dim2} have
looked at the algebraic structures that are required in order to be
able to define a exterior differential structure

\item[(\ref{dg_fn_genus})] 
Klimek and Lesniewski \cite{Klimek_Lesn2} have recently given an
outline of how one can create noncommutative analogues of
2-dimensional manifolds with genus greater than 2. They use the
noncommutative disc, so any new results for that should be useful.

\item[(\ref{dg_QG_top})]
It is believed that the true version of quantum gravity will naturally
lead to topology change.

\end{list}



\end{document}